\documentclass{mc}

\renewcommand\thisyear{202x}

\renewcommand\thisvolume{xx}
\renewcommand\datereceived{}
\renewcommand\dateaccepted{??}
\newcommand{\mathd}{\mathrm{d}}
\usepackage[numbers]{natbib}
\usepackage{mathrsfs}

\begin{document}

\markboth{Stephen Crowley}{Gaussian Processes Generated By Monotonically Modulated Stationary Kernels}

\title[Modulated Gaussian Processes]{Gaussian Processes Generated By Monotonically Modulated Stationary Kernels}

\author[Stephen Crowley]{Stephen Crowley\corrauth}
\email{{\tt stephencrowley214@gmail.com} (Stephen Crowley)}

\begin{abstract}
This article examines Gaussian processes generated by monotonically modulating stationary kernels. An explicit isometry 
between the original and the modulated reproducing kernel Hilbert spaces is established, preserving eigenvalues and normalization. 
The expected number of zeros over the interval $[0,T]$ is shown to be exactly $\sqrt{-\ddot{K}(0)}(\theta(T)-\theta(0))$, 
where $\ddot{K}(0)$ is the second derivative of the kernel at zero and $\theta$ is the modulating function.
\end{abstract}



\maketitle

\section{Introduction}

This article explores the properties of Gaussian
processes\cite{correlationTheoryOfStationaryRandomProcesses}\cite{stationaryAndRelatedStochasticProcesses}
generated by
monotonically modulating the kernels of stationary Gaussian processes. The
investigation centers on three key aspects: (1) the relationship between
eigenfunctions of the covariance operators defined by the original and the modulated kernels, 
(2) the preservation of normalization and eigenvalues under modulation, and (3) 
the expected number of zeros of the resulting processes. Beginning with a 
precise definition of the class of modulating functions $\mathcal{F}$, the 
article proceeds to establish theorems on eigenfunction transformation, normalization
preservation, and a formula for the expected value of the zero-counting function over $[0,T]$. 
These results provide a foundatio for understanding how stationary Gaussian processes
transform when modulated by monotonically increasing functions.
\section{Main Results}

\begin{definition}
  \label{scalingFunctions}Let $\mathcal{F}$ denote the class of functions
  $\theta : \mathbb{R} \to \mathbb{R}$ which are:
  \begin{enumerate}
    \item piecewise continuous with piecewise continuous first derivative,
    
    \item strictly monotonically increasing
    \begin{equation}
      \theta (t) < \theta (s) \forall - \infty \leqslant t < s \leqslant
      \infty
    \end{equation}
    \item and have a finite limiting derivative at infinity
    \begin{equation}
      \lim_{t \to \infty}  \dot{\theta} (t) < \infty
    \end{equation}
  \end{enumerate}
\end{definition}

\begin{remark}
  The conditions in Definition \ref{scalingFunctions} are somewhat redundant
  since a strictly monotonically increasing function must necessarily have a
  positive derivative.
\end{remark}

\

\begin{theorem}[Eigenfunctions]
  For any stationary kernel $K (t, s) = K (|t - s|)$, the eigenfunctions of
  the integral covariance operator
  \begin{equation}
    T_{K_{\theta}}  [f] (t) = \int_0^{\infty} K_{\theta} (| t - s |) f (s)
    \mathd s
  \end{equation}
  defined by the $\theta$-modulated kernel
  \begin{equation}
    K_{\theta} (t, s) = K (| \theta (t) - \theta (s) |)
  \end{equation}
  are given $\forall \theta \in \mathcal{F}$ by
  \begin{equation}
    \begin{array}{llllll}
      \phi_n (t) = \psi_n (\theta (t)) \sqrt{\dot{\theta} (t)} &  &  &  &  & 
    \end{array}
  \end{equation}
  which satisfies the eigenfunction equation
  \begin{equation}
    \begin{array}{ll}
      T_{K_{\theta}} [\phi_n] (t) & = \lambda_n \int_0^{\infty} K_{\theta} (|
      t - s |) \phi_n (s) \mathd s\\
      & = \lambda_n \int_0^{\infty} K_{\theta} (| t - s |) \psi_n (\theta
      (s)) \sqrt{\dot{\theta} (s)} \mathd s\\
      & = \lambda_n \int_0^{\infty} K (| \theta (t) - \theta (s) |) \psi_n
      (\theta (s)) \sqrt{\dot{\theta} (s)} \mathd s\\
      & = \lambda_n \phi_n (t)
    \end{array}
  \end{equation}
  where $\psi_n$ are the normalized eigenfunctions of the covariance operator defined by the original unmodulated
  kernel $K (|t - s|)$ which satisfy
  \begin{equation}
    \begin{array}{ll}
      T_K [\psi_n] (t) & = \lambda_n \int_0^{\infty} K  (| t - s |) \psi_n (s)
      \mathd s\\
      & = \lambda_n \psi_n (t)
    \end{array}
  \end{equation}
\end{theorem}

\begin{proof}
  The eigenfunction equation for the modulated kernel's covariance operator is:
  \begin{equation}
    \int_{- \infty}^{\infty} K (| \theta (t) - \theta (s) |) \phi_n (s) ds =
    \lambda_n \phi_n (t)
  \end{equation}
  The variables can be changed by substituting $u = \theta (s)$, $v = \theta
  (t)$:
  \begin{equation}
    \int_{- \infty}^{\infty} K (|v - u|) \frac{\phi_n (\theta^{- 1}
    (u))}{\dot{\theta} (\theta^{- 1} (u))} du = \lambda_n \phi_n (\theta^{- 1}
    (v))
  \end{equation}
  which is valid due to the strict monotonicity of $\theta$ which assures its
  invertability. Let
  \begin{equation}
    \psi_n (u) = \frac{\phi_n (\theta^{- 1} (u))}{\sqrt{\dot{\theta}
    (\theta^{- 1} (u))}}
  \end{equation}
  Then:
  \begin{equation}
    \int_{- \infty}^{\infty} K (|v - u|) \psi_n (u) du = \lambda_n \psi_n (v)
  \end{equation}
  This is precisely the eigenfunction equation for the original kernel $K$'s covariance operator. Therefore,
  \begin{equation}
    \phi_n (t) = \psi_n (\theta (t)) \sqrt{\dot{\theta} (t)}
  \end{equation}
  are the eigenfunctions of the modulated kernel's covariance operator 
  \begin{equation}
    \begin{array}{ll}
      T_{K_{\theta}} [\phi_n] (t) & = \lambda_n \int_0^{\infty} K_{\theta} (|
      t - s |) \phi_n (s) \mathd s
    \end{array}
  \end{equation}
  and $\psi_n$ are the eigenfunctions of the original kernel's covariance operator which satisfy
  \begin{equation}
    \begin{array}{ll}
      T_K [\psi_n] (t) & = \lambda_n \int_0^{\infty} K  (| t - s |) \psi_n (s)
      \mathd s
    \end{array}
  \end{equation}
  
\end{proof}

\begin{corollary}[Eigenvalue Invariance]
  The eigenvalues $\{\lambda_n \}$ of the modulated kernel $K_{\theta}$'s covariance operator are identical to those of the original kernel $K$'s covariance operator.
\end{corollary}

\begin{proof}
  For normalized $\psi_n$:
  \begin{equation}
    \int_{- \infty}^{\infty} | \phi_n (t) |^2 dt = \int_{- \infty}^{\infty} |
    \psi_n (\theta (t)) |^2 \dot{\theta} (t) dt
  \end{equation}
  Under the change of variables $u = \theta (t)$:
  \begin{equation}
    \int_{- \infty}^{\infty} | \psi_n (u) |^2 du = 1
  \end{equation}
  Therefore the $\phi_n$ are already normalized without additional constants.
\end{proof}

\begin{theorem}[Operator Conjugation]
  The transformation operator
  \begin{equation}
    M_{\theta} [\phi] (t) = \sqrt{\dot{\theta} (t)} \phi (\theta (t))
  \end{equation}
  conjugates the integral covariance operator
  \begin{equation}
    T_K [\phi] (t) = \int_0^{\infty} K (|t - s|) \phi (s) \mathd s
  \end{equation}
  where the resulting conjugated operator is
  \begin{equation}
    \begin{array}{ll}
      T_{K_{\theta}}  [\phi] (t) & {= M_{\theta}}  [T_K [M_{\theta}^{- 1}
      [\phi]]] (t)\\
      & = M \left[ \int_0^{\infty} K (|t - s|) \frac{\phi (\theta^{- 1}
      (s))}{\sqrt{\dot{\theta} (\theta^{- 1} (s))}} \mathd s \right] (t)\\
      & = \sqrt{\dot{\theta} (t)}  \int_0^{\infty} K (| \theta (t) - \theta (s)|)
      \frac{\phi (\theta^{- 1} (s))}{\sqrt{\dot{\theta} (\theta^{- 1} (s))}}
      \mathd s\\
      & = \int_0^{\infty} K (| \theta (t) - \theta (s) |) \phi (s) \mathd s\\
      & = \int_0^{\infty} K_{\theta} (| t - s |) \phi (s) \mathd s
    \end{array} \label{a}
  \end{equation}
  providing an explicit isometry between the original and modulated kernel
  Hilbert spaces.
\end{theorem}

\begin{proof}
  Observe that $M$ has inverse operator
  \begin{equation}
    M^{- 1} [\phi] (t) = \frac{\phi (\theta^{- 1} (t))}{\sqrt{\dot{\theta}
    (\theta^{- 1} (t))}}
  \end{equation}
  which follows from the invertibility of $\theta$ due to strict monotonicity
  and note that the last equality in Equation (\ref{a}) follows from the
  change of variables $s \mapsto \theta (s)$ with Jacobian $\dot{\theta} (s)$,
  demonstrating that the conjugated operator is precisely the integral
  operator with modulated kernel $K (| \theta (t) - \theta (s) |)$.
\end{proof}

\begin{theorem}[Expected Zero-Counting Function]
  Let $\theta \in \mathcal{F}$ and let $K (\cdot)$ be any positive-definite,
  stationary covariance function, twice differentiable at $0$. Consider the
  centered Gaussian process with covariance
  \begin{equation}
    K_{\theta} (t, s) = K (| \theta (t) - \theta (s) |)
  \end{equation}
  Then the expected number of zeros in $[0, T]$ is
  \begin{equation}
    \mathbb{E} [N ([0, T])] = \sqrt{- \ddot{K} (0)}  \hspace{0.17em} (\theta
    (T) - \theta (0))
  \end{equation}
\end{theorem}

\begin{proof}
  By the Kac-Rice
  formula{\cite[10.3.1]{stationaryAndRelatedStochasticProcesses}}:
  \begin{equation}
    \mathbb{E} [N ([0, T])] = \int_0^T \sqrt{- \lim_{s \to t}  \hspace{0.17em}
    \frac{\partial^2}{\partial t \partial s}  \hspace{0.17em} K_{\theta} (s,
    t)}  \hspace{0.27em} dt
  \end{equation}
  Computing the mixed partial derivative and taking the limit as $s \to t$:
  \begin{equation}
    \lim_{s \to t}  \hspace{0.17em} \frac{\partial^2}{\partial t \partial s} 
    \hspace{0.17em} K_{\theta} (s, t) = - \ddot{K} (0) \hspace{0.17em}
    \dot{\theta} (t)^2
  \end{equation}
  Therefore
  \begin{equation}
    \mathbb{E} [N ([0, T])] = \sqrt{- \ddot{K} (0)}  \int_0^T \dot{\theta} (t)
    \hspace{0.17em} dt = \sqrt{- \ddot{K} (0)}  \hspace{0.17em} (\theta (T) -
    \theta (0))
  \end{equation}
  so that
  \begin{equation}
    \begin{array}{ll}
      \sqrt{- \ddot{K} (0)}  \hspace{0.17em} (\theta (T) - \theta (0)) & =
      \sqrt{- \ddot{K} (0)}  \int_0^T \dot{\theta} (t)  \hspace{0.17em} dt\\
      & = \int_0^T \sqrt{- \ddot{K} (0) \dot{\theta} (t)^2}  \hspace{0.27em}
      dt\\
      & = \int_0^T \sqrt{- \lim_{s \to t}  \hspace{0.17em}
      \frac{\partial^2}{\partial t \partial s}  \hspace{0.17em} K (| \theta
      (t) - \theta (s) |)}  \hspace{0.27em} dt
    \end{array}
  \end{equation}
  which is precisely the Kac-Rice formula for the expected zero-counting
  function.
\end{proof}

\section{Conclusion}

The analysis presented in this article establishes several fundamental
properties of Gaussian processes generated by monotonically modulated
stationary kernels. Key results include: (1) a theorem demonstrating that the
eigenfunctions of the covariance operator defined by the modulated kernel are 
compositions of the stationary kernel's covariance operator eigenfunctions with 
the modulating function, times the square root of the modulating function's 
derivative, (2) proof of normalization and eigenvalue preservation under this 
transformation, establishing an isometry between the original and the modulated 
reproducing kernel Hilbert spaces, and (3) a concise formula for the expected value 
of the zero-counting function of the monotonically transformed process, expressed 
in terms of the original kernel's second derivative at zero times the modulating function's
values at the boundaries of the interval to which the expectation corresponds.


\bibliography{refs}

\end{document}